\documentclass[reqno,11pt]{article}
\usepackage{amsmath,amssymb}
\usepackage{geometry}
\usepackage{amsmath, amsfonts,amsthm,amssymb,amsbsy,upref,color,graphicx,amscd,enumerate}
\usepackage[active]{srcltx}
\usepackage[utf8]{inputenc}
\usepackage{algorithm}
\usepackage{stmaryrd}
\usepackage{array}
\usepackage{algorithmicx}
\usepackage{algpseudocode}
\usepackage{graphicx}
\usepackage{color}
\usepackage[table,xcdraw,svgnames]{xcolor}
\definecolor{light-salmon}{RGB}{255,140,120}
\graphicspath{{./metapost/}{./pics/}}

\usepackage[colorlinks]{hyperref}
\hypersetup{
	allbordercolors=.,
	citecolor=DarkGreen,
	linkcolor=DarkRed,
	urlcolor=NavyBlue
}

\usepackage{enumitem}

\theoremstyle{plain}
\newtheorem{thm}{Theorem}

\newtheorem{prop}[thm]{Proposition}

\theoremstyle{definition}
\newtheorem{rem}[thm]{Remark}

\newcommand{\Per}{\operatorname{Per}}

\renewcommand{\Bbb}{\mathbb}
\newcommand{\bo}[1]{{\bf #1}}

\makeatletter
\DeclareFontFamily{U}{tipa}{}
\DeclareFontShape{U}{tipa}{m}{n}{<->tipa10}{}
\newcommand{\arc@char}{{\usefont{U}{tipa}{m}{n}\symbol{62}}}%

\newcommand{\arc}[1]{\mathpalette\arc@arc{#1}}

\newcommand{\arc@arc}[2]{%
	\sbox0{$\m@th#1#2$}%
	\vbox{
		\hbox{\resizebox{\wd0}{\height}{\arc@char}}
		\nointerlineskip
		\box0
	}%
}
\makeatother

\title{A reverse isoperimetric inequality for convex shapes with inclusion constraint}

\author{Beniamin Bogosel}

\begin{document}
	 \maketitle
	 
	 \begin{abstract}
	 The convex shape contained in a disk having prescribed area and maximal perimeter is completely characterized in terms of the area fraction. The solution is always a polygon having all but one sides equal. The lengths of the sides are characterized through explicit equations. The case of more general containing shapes is also discussed from both theoretical and numerical perspectives.
	 \end{abstract}
	 
\section{Introduction}

The classical isoperimetric inequality in the plane states that the disk is the unique shape having fixed area and minimizing the perimeter. A survey of this problem can be found in \cite{Isop-AMM}. The reverse problem, maximizing the perimeter of shapes having fixed area, does not have solutions without adding extra constraints. Indeed, a sequence of flattening rectangles with sides $(a,1/a)$ does not have an upper bound for the perimeter. Nevertheless, if the shapes have prescribed area and are contained in a bounded set then a solution exists and its study is described below. 

Reverse isoperimetric inequalities are studied in various contexts. In \cite{Ball1991} it is shown that for any shape $\Omega$ there exists an affinely equivalent shape $\tilde \Omega$ such that the isoperimetric ratio $\Per(\tilde \Omega)/|\tilde \Omega|^{(n-1)/n}$ is smaller than the one given by the regular tetrahedron. The maximization of the perimeter at fixed area for shapes with bounds on the curvature are given in \cite{Howard1995} and \cite{Croce-reverse}. Investigating second order optimality conditions the authors of \cite{Lamboley2010} show that solutions to various shape optimization problems are polygonal. The case of the reverse isoperimetric problem naturally enters into this framework.

Considering convex shapes with given area contained in a disk, a complete description of perimeter maximizers can be achieved. Denote by $D$ the unit disk in the plane. For convex shapes $\Omega \subset D$ consider the problem 
\begin{equation}\label{eq:max-per}
 \max_{ |\Omega| =A} \Per(\Omega),
 \end{equation}
where $A \in (0,\pi)$ is given. The Blaschke selection theorem \cite[Theorem 1.8.7]{convex_bodies_Schneider} and the stability of the Hausdorff convergence of open sets for the inclusion \cite[p.33]{henrot-pierre-english} implies that \eqref{eq:max-per} has a solution $\Omega^*$. In the following it is proved that $\Omega^*$ is a polygon inscribed in $D$ which is completely characterized, up to a permutation of the sides, by the value of the area constraint. A similar result is mentioned in \cite{Favard}, however the proofs given below are slightly different and completely rigorous. In particular, Favard shows that the optimal shape should not have an arc in common with $\partial D$. This does not imply, however, that the optimal shape is polygonal. Situations where the optimal shape has a sequence of vertices on $\partial D$ with an accumulation point should be excluded.

The Lagrangian formulation of \eqref{eq:max-per} is described in \cite{bianchini_henrot} where problems of the form
\begin{equation}\label{eq:lagrange} \min_\Omega (\lambda|\Omega|-\Per(\Omega)),\end{equation}
are studied with $\Omega$ contained in an annulus. However, as underlined in \cite[Section 5]{bianchini_henrot} problems \eqref{eq:max-per} and \eqref{eq:lagrange} are not equivalent and techniques from \cite{bianchini_henrot} do not apply directly. The case where the shapes have prescribed area and a prescribed circumradius is described in \cite{Favard}. Although it is possible to use ideas from \cite{Favard} to deduce the result, a different strategy is employed below. First, the problem is considered in the class of $n$-gons with an upper bound on the number of sides. It is proved that the solution does not change, up to a permutation of the sides, when the number of sides becomes large enough. A classical polygonal approximation argument implies the result for problem \eqref{eq:max-per}.

First, let us restrict to $n$-gons. For $n \geq 3$ denote by $\mathcal P_n(A)$ the class of simple, convex polygons contained in $D$ with at most $n$ vertices, having area less than or equal to $A \in (0,\pi)$. It is obvious that $\mathcal P_n(A)$ is a compact set, described by a finite number of bounded parameters,  therefore
	 \begin{equation}\label{eq:max-per-poly}
	 \max_{ \Omega \in \mathcal P_n(A)} \Per(\Omega),
	 \end{equation}
also has a solution. The key to understanding \eqref{eq:max-per} is observing that \eqref{eq:max-per-poly} has a solution which does not change for $n \geq n_0$.

First let us investigate the class of admissible polygons $\mathcal P_n(A)$. If $P$ is a general $n$-gon which has a free vertex inside $D$ it is clear that moving it away from the adjacent diagonal towards $\partial D$ increases its area. Thus $n$-gons contained in $D$ with maximal area are inscribed $n$-gons. Considering central angles $\theta_i$, maximizing the area amounts to maximizing $\sum_{i=1}^n \sin \theta_i$ under the constraint $\sum_{i=1}^n \theta_i = 2\pi$. The concavity of the sine function and Jensen's inequality imply that the inscribed regular $n$-gon has maximal area among $n$-gons contained in $D$. Thus $\mathcal P_n(A)$ is non-void if and only if $A \leq \frac{n}{2}\sin \frac{2\pi}{n}$.

Let $P^*$ be a solution for \eqref{eq:max-per-poly} and suppose that $P^*$ has $k\leq n$ vertices. Without loss of generality, suppose that no vertex in $P^*$ is redundant, i.e. the vertices are distinct and no three vertices are colinear. Since $|P^*|\leq A>0$, $P^*$ has at least three vertices. Indeed, if $|P^*|=0$ then $P^*$ would be reduced to a segment. Replacing it with a thin symmetric quadrilateral close to the diameter still fits in $D$ and obviously increases its perimeter, while remaining admissible.

The following result is immediate:

\begin{thm}\label{thm:free-vert}
	Let $P^*$ be a solution of \eqref{eq:max-per-poly}. Then $P^*$ is inscribed in $D$, i.e., there are no \emph{free vertices} among vertices of $P^*$.
\end{thm} 

\begin{figure}
	\centering 
	\includegraphics[width=0.4\textwidth]{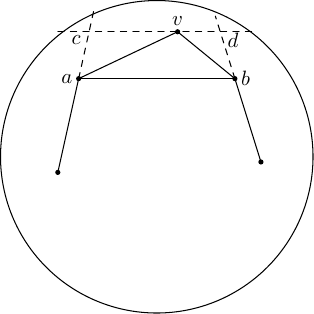}
	\caption{An optimal $n$-gon does not have a free vertex inside the container $D$.}
	\label{fig:freevert}
\end{figure}

\emph{Proof:} Take $v$ a vertex of $P^*$ and suppose $v \notin \partial D$. Consider $a,b$ the vertices adjacent to $v$ in $P^*$. The line parallel to $ab$ through $v$ generates a chord in $D$. Denote $[cd]$ the intersection of this chord with the region generated by the edges of $P^*$ adjacent to $av, vb$. See Figure \ref{fig:freevert} for an illustration. Since the function 
\[ [cd] \ni v \mapsto |va|+|vb|\]
is strictly convex, its maximum is attained for extremal points. Thus, either $v$ belongs to $\partial D$ or $va, vb$ are colinear with one of the adjacent sides, showing that $a$ or $b$ are redundant, contradicting the assumptions. \hfill $\square$

The previous process shows that no convex polygon in $\mathcal P_n$ which has a vertex interior to $D$ is optimal for \eqref{eq:max-per-poly}. Thus, a solution $P^*$ for \eqref{eq:max-per-poly} is an inscribed polygon in $D$. The same procedure applies to general containing shapes leading to the following result:

\begin{thm}\label{thm:optimum-inscribed}
	Let $\Omega$ be a bounded convex shape, not reduced to a segment. Let $P^*$ be a polygon maximizing the perimeter among polygons $P$ with at most $n$ sides are contained in $\Omega$ and having fixed area $|P|=A \in (0,|\Omega|)$. Then $P^*$ is inscribed in $\Omega$, i.e., all vertices of $P^*$ belong to $\partial \Omega$.
\end{thm}

Moreover, in the case of the disk, let us show that $P^*$ must contain the center of $D$ in its interior. If $P^*$ does not contain the center of $D$ then $P^*$ is contained in the interior of a half-disk and may be translated such that $P^*$ is strictly in the interior of $D$. Then at least one of the vertices of $P^*$ can be moved parallel to the adjacent diagonal contradicting again the optimality of $P^*$.

An inscribed polygon $P^*$ in the disk $D$ is completely determined, up to a permutation, by the lengths of its sides, or equivalently, the angles at the center made by the sides. 

Take $3\leq m \leq n$ and consider $\theta_1,...,\theta_m$ the corresponding center angles. Since the center of $D$ is interior to the polygon, we have $\theta_i \in [0,\pi]$, $i=1,...,m$. The area of the inscribed polygon is
\[ |P| = \frac{1}{2} \sum_{i=1}^m \sin \theta_i\]
and its perimeter is 
\[ \Per(P) = 2\sum_{i=1}^m \sin \frac{\theta_i}{2}.\]
Of course, the central angles also verify the relation $\sum_{i=1}^m \theta_i = 2\pi$. 

The optimal polygon $P^*$ must saturate the area constraint. While heuristic arguments can be given, this is a direct consequence of Proposition \ref{prop:general-ineq}, part (c). Thus, angles characterizing $P^*$ are solutions to the problem 
\begin{equation}\label{eq:max-discrete}
\max_{(\theta_i) \in X_m} \sum_{i=1}^m 2\sin \frac{\theta_i}{2}
\end{equation}
where 
\[X_m = \{ (\theta_i)_{i=1}^m : \sum_{i=1}^m \theta_i = 2\pi, \sum_{i=1}^m \sin \theta_i = 2A\}.\]

In Section \ref{sec:main-results} a series of elementary results in analysis and constrained optimization help characterize the solutions of problem \eqref{eq:max-discrete}. The case of a general container is investigated theoretically in Section \ref{sec:general} and numerical simulations are presented in Section \ref{sec:numerics}.





\section{Main results}
\label{sec:main-results}

\subsection{Analysis of a family of constrained optimization problems}

In this section some elementary optimization problems are studied which will help prove the main results of the paper.

\begin{prop}\label{prop:ineq-3}
	(a) Given $S \in [0,2\pi]$ and $T\in (0,3\sin \frac{S}{3}]$, suppose that $0 < a \leq b <c\leq \pi$ are given, such that 
	\begin{equation}\label{eq:constraint-3}
	 a+b+c = S,\ \sin a+\sin b +\sin c = T.
	 \end{equation}
	
	 Then there exists $\varepsilon>0$ such that for every $s \in [0,\varepsilon]$ there exists $0<t(s)<\frac{c-b-s}{2}$ such that
	\begin{equation}\label{eq:equation-t}
	 \sin (a-s) + \sin(b+s+t(s))+\sin (c-t(s))= T.
	 \end{equation}
	
	(b) If $0 \leq a\leq b\leq c\leq \pi$ verify the constraints \eqref{eq:constraint-3} and they maximize 
	\[ \sin \frac{a}{2}+\sin \frac{b}{2}+\sin \frac{c}{2}\]
	then either $a=0$ or $b=c$.
\end{prop}

The condition $T\leq 3\sin \frac{S}{3}$ is necessary since $x \mapsto \sin x$ is concave on $[0,\pi]$ and the upper bound for $\sin a +\sin b +\sin c$ follows using Jensen's inequality.

\emph{Proof:} (a) We have $a,b,c \in [0,\pi]$ and the sine function is concave on this interval. It is straightforward to see that mappings of the form 
\[ \tau \mapsto \sin(x+\tau)+\sin(y-\tau), 0<x<y\leq 2\pi \]
are strictly increasing. Simply differentiate and use the concavity of the sine function. Therefore, for $s>0$ small enough we have
\[ \sin(a-s)+\sin(b+s)<\sin a+\sin b\]
and
\[q: t \mapsto  \sin(a-s)+\sin(b+s+t)+\sin(c-t)\]
is increasing on $[0,\frac{c-b-s}{2}]$. Since $b<c$, the concavity of the sine function and Jensen's inequality implies that $\sin a+2\sin \frac{b+c}{2}>T$.

 Thus, there exists $\varepsilon>0$ and an interval $[0,\varepsilon]$ such that $q(\frac{c-b-s}{2}) = \sin(a-s)+2\sin \frac{b+c+s}{2}>T$. The continuity of $q$ implies the existence of $t(s)$ verifying \eqref{eq:equation-t}. 

(b) If $t(s)$ verifies (a) then the partial derivative of \eqref{eq:equation-t} with respect to $t$ is non-zero. The implicit function theorem implies that the mapping $s \mapsto t(s)$ is differentiable. Differentiating \eqref{eq:equation-t} with respect to $s$ gives
\begin{equation}\label{eq:deriv-t}
t'(s) = \frac{\cos(a-s)-\cos(b+s+t(s))}{\cos(b+s+t(s))-\cos(c-t(s))}.
\end{equation}
Of course, $t'(s)>0$ since $\cos$ is decreasing on $[0,\pi]$.  

Consider 
\begin{equation}\label{eq:def-h} h: s \mapsto \sin\frac{a-s}{2}+\sin\frac{b+s+t(s)}{2}+\sin \frac{c-t(s)}{2}.
\end{equation}
The derivative with respect to $s$ is
\[ h'(s) = \frac{1}{2} \left( \cos \frac{b+s+t(s)}{2}-\cos \frac{a-s}{2}\right) + \frac{1}{2}t'(s) \left(\cos \frac{b+s+t(s)}{2}-\cos \frac{c-t(s)}{2} \right).\]
Replacing $t'(s)$ from \eqref{eq:deriv-t} shows that $h'(s)>0$ is equivalent to 
\begin{equation}\label{eq:equiv-deriv-pos} \frac{\cos \frac{b+s+t(s)}{2}-\cos\frac{c-s}{2}  }{\cos(b+s+t(s))-\cos(c-t(s))}>\frac{\cos \frac{a-s}{2}-\cos \frac{b+s+t(s)}{2}}{\cos(a-s)-\cos(b+s+t(s))}.  
\end{equation}
Both ratios are of the form
\[ \frac{\cos \frac{x}{2}-\cos \frac{y}{2}}{\cos x-\cos y} = \frac{1}{2} \frac{1}{\cos x+\cos y},\]
where the well known equality $\cos x = 2\cos^2 \frac{x}{2}-1$ was used. Therefore \eqref{eq:equiv-deriv-pos} is equivalent to
\[ \frac{1}{\cos (b+s+t(s))+\cos(c-t(s))}> \frac{1}{\cos (a-s)+\cos(b+s+t(s)},\]
which is true, since $\cos$ is decreasing on $[0,\pi]$ and $0\leq a-s\leq b+s+t(s)< c-t(s)\leq \pi$. Thus $h$ is strictly increasing. 

Consider the problem of maximizing 
\[ \sin \frac{a}{2}+\sin \frac{b}{2}+\sin \frac{c}{2}\]
such that $a+b+c=S$, $a,b,c \in [0,\pi]$ and $\sin a+\sin b+\sin c = T$. Suppose that $0<T<S\leq 2\pi$ are given such that the set of admissible triples $(a,b,c)$ is non-void. The set of admissible points is non-void and compact and the function to be maximized is continuous. Therefore, there exists a solution $(a,b,c)$ ordered such that $a\leq b \leq c$.

If $a>0$ and $b<c$ then the results in (a) imply that $h$ given by \eqref{eq:def-h} is strictly increasing and a perturbation of the form
\[ 0 \leq a-s\leq b+s+t(s)\leq c-t(s)\leq \pi, s>0,t(s)>0\]
strictly increases the objective function. Therefore, a maximizer verifies $a=0$ or $b=c$.
\hfill $\square$

Following similar ideas to those in Proposition \ref{prop:ineq-3} the following, more general result can be proved.

\begin{prop}\label{prop:general}
	Consider functions $f,g:[0,M] \to \Bbb{R}_+$ of class $C^2$ with $f',g'$ injective and $g''\neq 0$ on $[0,M]$. Moreover, assume that $f''/g''$ is injective on $[0,M]$. Consider
	\[ X = \{ a,b,c \in [0,M]: a+b+c = S, g(a)+g(b)+g(c) = T\},\]
	where $S,T$ are chosen such that $X$ is not empty or reduced to a point. 
	
	Then if $(a,b,c)$ is a solution to the problem 
	\[ \min_{(a,b,c) \in X} f(a)+f(b)+f(c)\]
	such that $a,b,c \in (0,M)$ then at least two of $a,b,c$ are equal.
\end{prop}

\emph{Proof:} The optimization problem stated previously deals with the minimization of a continuous function on a compact set, therefore solutions exist. Let $(a,b,c) \in X$ be a solution and suppose that $0<a<b<c<M$, without loss of generality. Then the Jacobian of the mapping $F(a,b,c) = \begin{pmatrix}
a+b+c \\
g(a)+g(b)+g(c)
\end{pmatrix}$ has rank two, since $g'$ is injective. The implicit function theorem states that there exists a parametrization of $X$ around $(a,b,c)$, given by $(-\varepsilon,\varepsilon) \ni t \to (a(t),b(t),c(t))$ such that $a,b,c$ are continuously differentiable with respect to $t$ and $(a(t),b(t),c(t)) \in X$ for $t \in (-\varepsilon,\varepsilon)$. Without loss of generality assume that $a'(t)>0$ on $(-\varepsilon,\varepsilon)$. In the following, we sometimes drop the variable $t$, for simplifying the notations. Differentiating the constraints with respect to $t$ gives
\[ a'+b'+c' = 0, g'(a)a'+g'(b)b'+g'(c)c' = 0.\]
Using $b' = -a'-c'$ gives
\begin{equation}\label{eq:constraint-g}
 c'(g'(c)-g'(b))-a'(g'(b)-g'(a)) = 0.
 \end{equation}

A classical result known as the Cauchy mean value theorem states that if $h_1,h_2$ are of class $C^1$, $h_2'\neq 0$ on an interval $[x,y]$ then there exists $z \in (x,y)$ such that
\[ \frac{h_1(y)-h_1(x)}{h_2(y)-h_2(x)} = \frac{h_1'(z)}{h_2'(z)}.\]
Apply this result for $h_1=f'$ and $h_2=g'$ like in the hypothesis on the intervals $[a,b],[b,c]$, obtaining $z_{ab} \in (a,b)$ and $z_{bc} \in (b,c)$ such that
\[ \frac{f'(b)-f'(a)}{g'(b)-g'(a)} = \frac{f''(z_{ab})}{g''(z_{ab})},\ \frac{f'(c)-f'(b)}{g'(c)-g'(b)} = \frac{f''(z_{bc})}{g''(z_{bc})}.\]

This allows us to perform the following computation:
\begin{align*}
f'(a)a'+f'(b)b'+f'(c)c' & = c'(f'(c)-f'(b))-a'(f'(b)-f'(a)) \\
& = a'\left[ \frac{c'}{a'}(f'(c)-f'(b))-(f'(b)-f'(a))\right] \\
& = a'(g'(b)-g'(a))\left[ \frac{f'(c)-f'(b)}{g'(c)-g'(b)} - \frac{f'(b)-f'(a)}{g'(b)-g'(a)}\right]\\
& = a'(g'(b)-g'(a))\left[\frac{f''(z_{bc})}{g''(z_{bc})}-\frac{f''(z_{ab})}{g''(z_{ab})} \right],
\end{align*}
where the constraints on $a,b,c$ and \eqref{eq:constraint-g} were used. The optimality of $a,b,c$ and $a'>0$ implies that either $a=b$ or $z_{bc}=z_{ab}$, both contradicting the fact that $a,b,c$ are distinct. \hfill $\square$

Previous results allow to state the following generalization of Proposition \ref{prop:ineq-3}.

\begin{prop}\label{prop:general-ineq}
	Let $n \geq 3$ and define 
	\[ X_{S,T,n}:=\left\{ (\theta_i)_{i=1}^n : 0\leq \theta_1\leq ... \leq \theta_n, \sum_{i=1}^n \theta_i = S, \sum_{i=1}^n \sin \theta_i = T\right\},\]
	where $S \in (0,2\pi],T\in (0,n\sin \frac{S}{n}]$ are given such that $X_{S,T,n}$ is non-void. 
	
	(a) Suppose $(\theta_i)_{i=1}^n$ is a solution of 
	\begin{equation}\label{eq:pb-max-n} \max_{(\theta_i)\in X_{S,T,n}} \sum_{i=1}^n \sin \frac{\theta_i}{2}.\end{equation}
	Then either $\theta_1=0$ or $0<\theta_1\leq \theta_2=\theta_3=...=\theta_n$.
	
	(b) Given $0<T<S\leq 2\pi$ such that $X_{S,T,n}$ is non-void, there exists $n_0$ such that the solution of \eqref{eq:pb-max-n} does not change (up to a permutation) for $n \geq n_0$.
	
	(c) The maximal value attained in \eqref{eq:pb-max-n} is strictly increasing with respect to $T \in (0,n\sin \frac{S}{n}]$.
\end{prop}

\emph{Proof:} (a) Suppose that $0<T<S\leq 2\pi$ are given such that $X_{S,T,n} \neq \emptyset$. It is obvious that $X_{S,T,n}$ is compact, therefore \eqref{eq:pb-max-n} has a solution.

Take three consecutive variables $(\theta_i,\theta_{i+1},\theta_{i+2})$ of a solution. In view of the results in Proposition \ref{prop:ineq-3} (b), either $\theta_i=0$ or $\theta_{i+1}=\theta_{i+2}$. 

If $\theta_i=0$ for some $i$ then $\theta_1=0$. Otherwise, $\theta_{i+1}=\theta_{i+2}$ for all $i \geq 1$ implying that $0<\theta_1 \leq \theta_2=...=\theta_n$.

(b) Fix $0<T<S\leq 2\pi$ such that $X_{S,T,n}\neq \emptyset$ for some $n\geq 3$. Note that adding some extra variables equal to zero does not change the constraint set or the objective function. Thus $X_{S,T,n}$ can be identified with a subset of $X_{S,T,m}$ for $m\geq n$.

If the solution of \eqref{eq:pb-max-n} verifies $0<\theta_1\leq \theta_2=...=\theta_n$ then the following system admits at least one solution:
\[ \theta_1+(n-1)\theta_2=S, \ \sin \theta_1+(n-1)\sin \theta_2 = T, 0<\theta_1\leq \theta_2.\]
Then $\theta_1 \in (0,\frac{S}{n}]$, $\theta_2 = \frac{S-\theta_1}{n-1}$ and $\theta_1$ verifies
\begin{equation}\label{eq:theta1} \sin \theta_1+(n-1)\sin \frac{S-\theta_1}{n-1} = T.
\end{equation}
The function 
\[ \theta \mapsto \sin \theta+(n-1)\sin \frac{S-\theta}{n-1},\]
is strictly increasing on $[0,S/n]$. Therefore, \eqref{eq:theta1} has a solution in $(0,S/n]$ if and only if 
\begin{equation}\label{eq:ineq-B}
 T \in \Big((n-1)\sin \frac{S}{n-1},n \sin \frac{S}{n}\Big].
\end{equation}
Note that $(n-1)\sin \frac{S}{n-1}<n \sin \frac{S}{n}$ since $0<S \leq 2\pi$, $n \geq 3$ and $x \mapsto \frac{\sin x}{x}$ is decreasing on $[0,\pi]$. Moreover, the sequence $n \sin \frac{S}{n}$ is increasing and converges to $S$.

Therefore, given $S \in (0,2\pi]$ and $T\in (0,S)$, there exists a unique $n$ such that \eqref{eq:theta1} has a solution in $(0,S/n]$ and that $n$ verifies \eqref{eq:ineq-B}. Therefore, there exists a unique $n_0$ such that \eqref{eq:pb-max-n} has a solution in $X_{S,T,n_0}$ with $\theta_1>0$. 

If $n\geq n_0$ then strictly positive components of the solution of \eqref{eq:pb-max-n} in $X_{S,T,n}$ are the same (up to a permutation) as the the ones in the solution in $X_{S,T,n_0}$. 

(c) The previous arguments show that $\theta_1$, solution of \eqref{eq:theta1}, is strictly increasing with respect to $T$. The objective function has the same montonicity properties, thus the maximal value in \eqref{eq:pb-max-n} is strictly increasing with respect to $T$.
\hfill $\square$

\subsection{Reverse isoperimetric inequality for convex polygons contained in a disk}

We are now ready to state and prove the main result of the paper. This result is also stated in \cite{Favard}, however the proof argument stated there is incomplete. Favard proves that the convex shape maximizing the perimeter cannot contain an arc of the boundary of the bounding disk $\partial D$. This does not show that the optimizing shape is a polygon. In addition it should be proved that an inscribed polygonal shape with vertices on $\partial D$ and an accumulation point for the vertices is not optimal.

\begin{thm}\label{thm:max-per}
	Solutions to problem \eqref{eq:max-per-poly} are uniquely characterized by $n$ and $A$, up to a permutation of the sides. Given $A \in (0,\pi)$, there exists $m$ such that for $n \geq m$ the solution to problem \eqref{eq:max-per-poly} does not change, up to a permutation of the sides. 
	
	Given $A \in (0,\pi)$, solution to problem \eqref{eq:max-per} is a polygon inscribed in $D$, having $m$ sides with $m$ verifying 
	\begin{equation}\label{eq:ineq-Area} (m-1)\sin \frac{2\pi}{m-1}<2A\leq m \sin \frac{2\pi}{m}.
	\end{equation}
	
	The angles subtended by the sides have the form $0<\theta_1\leq \theta_2=...=\theta_m$ where $\theta_1\in (0,\frac{2\pi}{m}]$ solves 
	\[ \sin \theta_1+(m-1)\sin \frac{2\pi-\theta_1}{m-1}=2A.\] 
\end{thm}

\emph{Proof:} In view of the discussion in the introduction, solutions to \eqref{eq:max-per-poly} are inscribed polygons in $D$. Considering the central angles $\theta_i$ associated to each sides, then $(\theta_i)_{i=1}^m$ solve the optimization problem \eqref{eq:max-discrete}. Proposition \ref{prop:general-ineq} applied for $S = 2\pi, T=2A$ shows that solutions of \eqref{eq:max-discrete} verify 
\[ 0<\theta_1\leq \theta_2=...=\theta_m,\]
where $m$ is the unique positive integer verifying \eqref{eq:ineq-Area}. Thus, for $n\geq m$, the solution of \eqref{eq:max-per-poly} does not change, up to a permutation of the sides.

Consider a solution $\Omega^*$ to problem \eqref{eq:max-per} and let $Q_n$ be a sequence of convex polygons converging to $\Omega^*$. Since the maximizer $\Omega^*$ contains the center of $D$, as recalled in the introduction, assume that $Q_n$ contain the origin. Without loss of generality we may assume $Q_n \in \mathcal P_n(A)$. Indeed, consider for each $Q_n$ a maximal homothety factor $\lambda_n \in (0,1]$ such that $\lambda_nQ_n$ has area at most $A$ and is contained in $D$. If $\lambda_n<1$ replace $Q_n$ with $\lambda_n Q_n$. Since the limit $\Omega^*$ has area $A$ and is contained in $D$, we have $\lambda_n \to 1$ as $n \to \infty$. For each $n$ consider a solution $R_n \in \mathcal P_n(A)$ for \eqref{eq:max-per-poly}. Since $R_n$ is inscribed in $D$, its sides of positive length are given by the central angles $0<\theta_1\leq \theta_2=...=\theta_m$ and there exist a finite number of permutations of these $m$ sides, the sequence $R_n$ has a converging subsequence which has a limit $\Omega_m \in \mathcal P_m$. The optimality of $R_n$ implies $\Per(Q_n)\leq \Per(R_n)$ and passing to the limit we obtain
\[ \Per(\Omega^*) \leq \Per(\Omega_m),\]
where $\Omega_m$ is a polygon in $\mathcal P_m(A)$. Therefore, there exists a polygonal solution to \eqref{eq:max-per} for any choice of $A \in (0,\pi)$.

The equation verified by the angles can be found in the proof of Proposition \ref{prop:general-ineq}. \hfill $\square$

Solutions of problem \eqref{eq:lagrange}, the Lagrangian formulation for \eqref{eq:max-per}, can also be characterized, assuming $\Omega \subset D$.

\begin{thm}
	For every $\lambda>0$, problem 
	\begin{equation}\label{eq:min-lagrange} \min_{\Omega \subset D} \lambda |\Omega|-\Per(\Omega),\end{equation}
	where $\Omega$ is assumed convex, 
	has a solution. The solution is either the disk $D$, a segment which is a diameter of $D$ or a polygon completely characterized by $\lambda$.
\end{thm}

\emph{Proof:} First, let us note that for $\lambda>0$ the objective function is bounded and the Blaschke selection theorem implies the existence of a minimizer $\Omega^*$. If $\Omega^*=D$ there is nothing to prove.

Suppose that $|\Omega^*|<|D|$. Then $\Omega^*$ also solves \eqref{eq:max-per} for $A=|\Omega^*|$. Results of Theorem \ref{thm:max-per} imply that $\Omega^*$ is a polygon completely characterized by its area. \hfill $\square$

\begin{figure}
	\centering
	\includegraphics[width=0.6\textwidth]{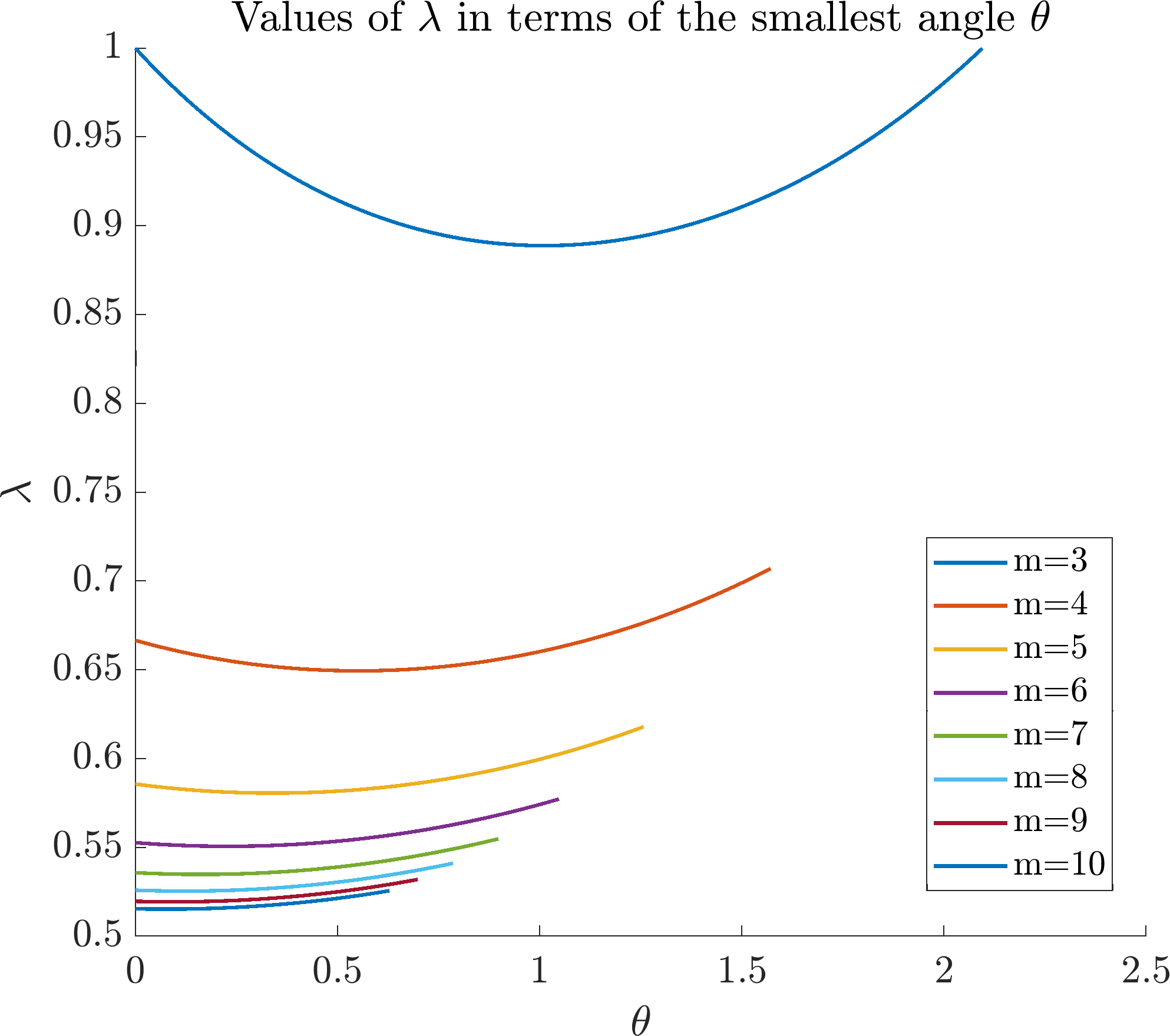}
	\caption{Values of $\lambda$ in terms of $m$ and $\theta$ given in \eqref{eq:lambdas} for $3 \leq m \leq 10$.}
	\label{fig:lambdas}
\end{figure}

We do not attempt to fully characterize solutions of \eqref{eq:min-lagrange} since this can be achieved using ideas in \cite{bianchini_henrot}. The solution $\Omega^*$ could be identified as follows. Consider $m$ the number of sides of $\Omega^*$, assuming it is a non-trivial $m$-gon. Consider the corresponding angles $0<\theta_1\leq \theta_2 =... = \theta_m$. Then, denoting $\theta = \theta_1\in (0,\frac{2\pi}{m}]$ we have
\[ \lambda |\Omega| - \Per(\Omega) = \lambda\frac{1}{2} \left( \sin \theta +(m-1) \sin \frac{2\pi-\theta}{m-1}\right) - 2\left(\sin \frac{\theta}{2}+(m-1)\sin \frac{2\pi-\theta}{2(m-1)}\right).\]
If $\Omega^*$ is not the regular $m$-gon, the optimality of $\Omega^*$ implies that $\theta$ is a critical point for the previous expression of $\theta$. Using $\cos x-\cos y = 2(\cos \frac{x}{2}-\cos \frac{y}{2})(\cos \frac{x}{2}+\cos \frac{y}{2})$ we obtain
\begin{equation}\label{eq:lambdas}
 \lambda = \left(\cos \frac{\theta}{2}+\cos \frac{2\pi-\theta}{2(m-1)}\right)^{-1},
 \end{equation}
where $\theta \in (0,\frac{2\pi}{m}]$. A plot of the values of $\lambda$ in terms of $m$ and $\theta$ is shown in Figure \ref{fig:lambdas} for $3 \leq m \leq 10$. It should be noted that not all values of $\lambda\geq 0$ can be written in the form given by \eqref{eq:lambdas}.

Numerical simulations suggest that the gaps between different branches of possible values for $\lambda$ correspond to regular polygons with the corresponding number of sides. If $\lambda>1$ the optimal shape is a diameter of the disk while for $\lambda\leq 0.5$ the solution is the whole disk. 

\section{General convex sets}
\label{sec:general}

The disk containing the convex shapes can be replaced by a general convex set. The existence of perimeter maximizers under area constraint is straightforward.

\begin{thm}
	Let $\Omega$ be a convex and bounded shape with non-empty interior and let $A \in (0,|\Omega|)$ be a choice for the area constraint. Then, assuming $\omega \subset \Omega$ is convex, the problem 
	\[ \max_{|\omega |=A,\ \omega \subset \Omega } \Per(\omega)\]
	has a solution. 
\end{thm}

\emph{Proof:} Like in the case of the disk, combine the Blaschke selection theorem \cite{convex_bodies_Schneider} and the stability of the Hausdorff convergence for set inclusion \cite[p.33]{henrot-pierre-english}. \hfill $\square$

\begin{rem}
	The convexity of the container $\Omega$ could be relaxed, since the inclusion constraint and convexity are preserved under Hausdorff convergence. Nevertheless, the study of the maximal perimeter convex subset is more difficult in this case since the optimal shape can touch non-convex parts in $\partial \Omega$ at most one point. 
\end{rem}

In the case when the container $\Omega$ is a convex polygon the following result holds.

\begin{thm}\label{thm:max-per-poly}
	Let $\Omega$ be a convex $n$-gon and $A \in (0,|\Omega|)$. Then there exist a solution of
		\begin{equation}\label{eq:polygon-container} \max_{|\omega |=A,\ \omega \subset \Omega } \Per(\omega)
		\end{equation}
		which is polygonal.
\end{thm} 

\emph{Proof:} Like in the case of the disk, consider $\mathcal P_n(A,\Omega)$ the family of convex polygons included in $\Omega$ with at most $m$ sides and area at most $A$. Consider the discrete problem
\begin{equation}\label{eq:discrete-poly-general}
 \max_{P \in \mathcal P_m(A,\Omega) } \Per(P).
 \end{equation}
This problem has solutions, since it is finite dimensional in nature and the perimeter depends continuously on the vertices of a polygon $P$. Like in the case of the disk (Theorem \ref{thm:free-vert}) an optimal polygon $P$ cannot have free vertices in the interior of $\Omega$. Therefore solutions of \eqref{eq:discrete-poly-general} are polygons $P$ whose vertices lie on the boundary $\partial \Omega$ of the containing set.

Since $\Omega$ is an $n$-gon, for each edge of $\Omega$ there are at most $2$ non-redundant vertices of an inscribed convex polygon $P$. Thus, solutions of \eqref{eq:discrete-poly-general} have at most $2n$ vertices. Therefore, for $m \geq 2n$ the solution $P_m$ to \eqref{eq:discrete-poly-general} is also a solution for $m=2n$. In any case, there exists a polygon $P^*$ with at most $2n$ sides which is solution to \eqref{eq:discrete-poly-general} for all $m \geq 2n$.

Now consider a solution $\omega$ for \eqref{eq:polygon-container}. Let $P_k$ be a polygonal approximation of $\omega$ such that $P_k \subset \omega$ for every $n$. For example, one might consider inscribed $n$-gons with uniformly distributed vertices on $\partial \omega$. Of course, we have $P_k\to \omega$ in the Hausdorff convergence when $k \to \infty$. For $k \geq 2n$ polygons $P_k$ are admissible in \eqref{eq:discrete-poly-general}, showing that $\Per(P_k) \leq \Per(P^*)$. Letting $k \to \infty$ we have $\Per(\omega) \leq \Per(P^*)$. Thus, $P^*$ is a polygonal solution for \eqref{eq:polygon-container}. \hfill $\square$

\begin{rem}
	It is tempting to conjecture that replacing $D$ with a convex or strictly convex set $\Omega$, we have the following result: convex sets $\omega \subset \Omega$ having prescribed area and maximal perimeter are polygons. An argument like in \cite{Favard} could be formalized: if the optimal set $\omega^*$ contains a segment $ab$ and a subset $\arc{bc}$ of $\partial \Omega$, then assuming $\arc{bc}$ is small enough, replacing $b$ with $b'$ such that the triangle $\Delta ab'c$ and the curvilinear triangle $abc$ have the same area increases the perimeter: $ab+\arc{bc}<ab'+b'c$. When $\Omega$ is not a disk it is difficult to investigate inscribed polygons with fixed area and maximal perimeter. If a parametrization of $\partial \Omega$ is available, the length and area bounded by a chord depends on the position of both endpoints, not only on the length of the chord, like in the case of the disk.
\end{rem}

\section{Numerical aspects}
\label{sec:numerics}

\subsection{The disk}

The characterization of the solutions of problem \eqref{eq:max-per} given in Theorem \ref{thm:max-per} suggests the following algorithm for finding the optimal shape, given a volume fraction $A \in (0,\pi)$. 

(i) Find the number of sides of the polygon by identifying $n \geq 2$ such that $ (n-1) \sin \frac{2\pi}{n-1} <2A \leq n \sin \frac{2\pi}{n}$. Since $q: x \mapsto \frac{\sin x}{x}$ is decreasing for $x \geq 0$, it is enough to find $x_A$ such that $q(x_A) = A/2\pi$, then use $\frac{2\pi}{n} \leq x_A <\frac{2\pi}{n-1}$ to find $n$.

(ii) Next, in view of \eqref{eq:theta1}, solve the equation 
\[  \sin \theta+(n-1)\sin \frac{2\pi-\theta}{n-1} = 2A\]
which has a unique solution $\theta$ in the interval $[0,\frac{2\pi}{n}]$. This gives the central angle of the smallest side. The other remaining edges all have equal lengths and central angles $(2\pi-\theta)/(n-1)$.

\begin{figure}
	\centering 
	\begin{tabular}{ccc}
		\includegraphics[width=0.3\textwidth]{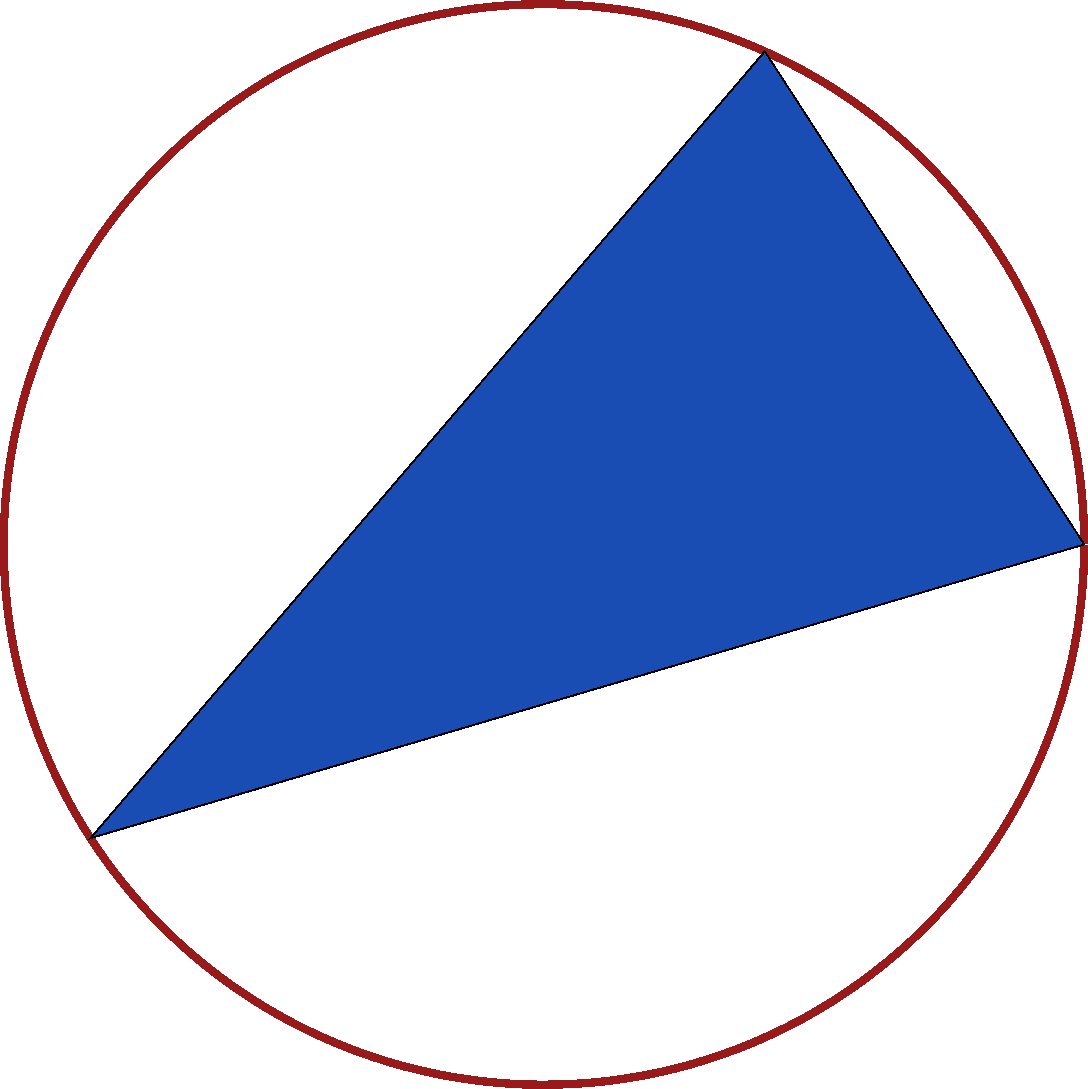}&
		\includegraphics[width=0.3\textwidth]{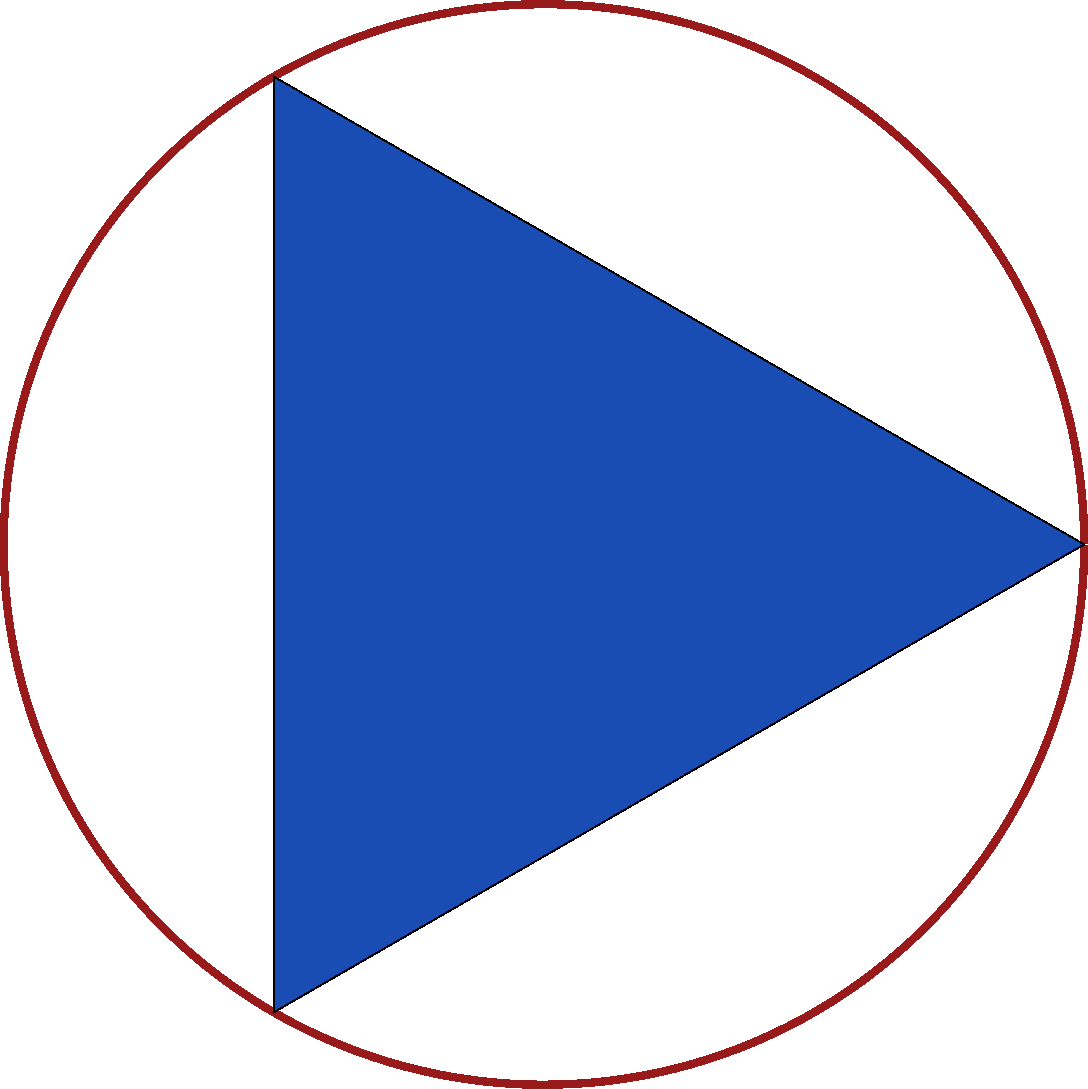}&
		\includegraphics[width=0.3\textwidth]{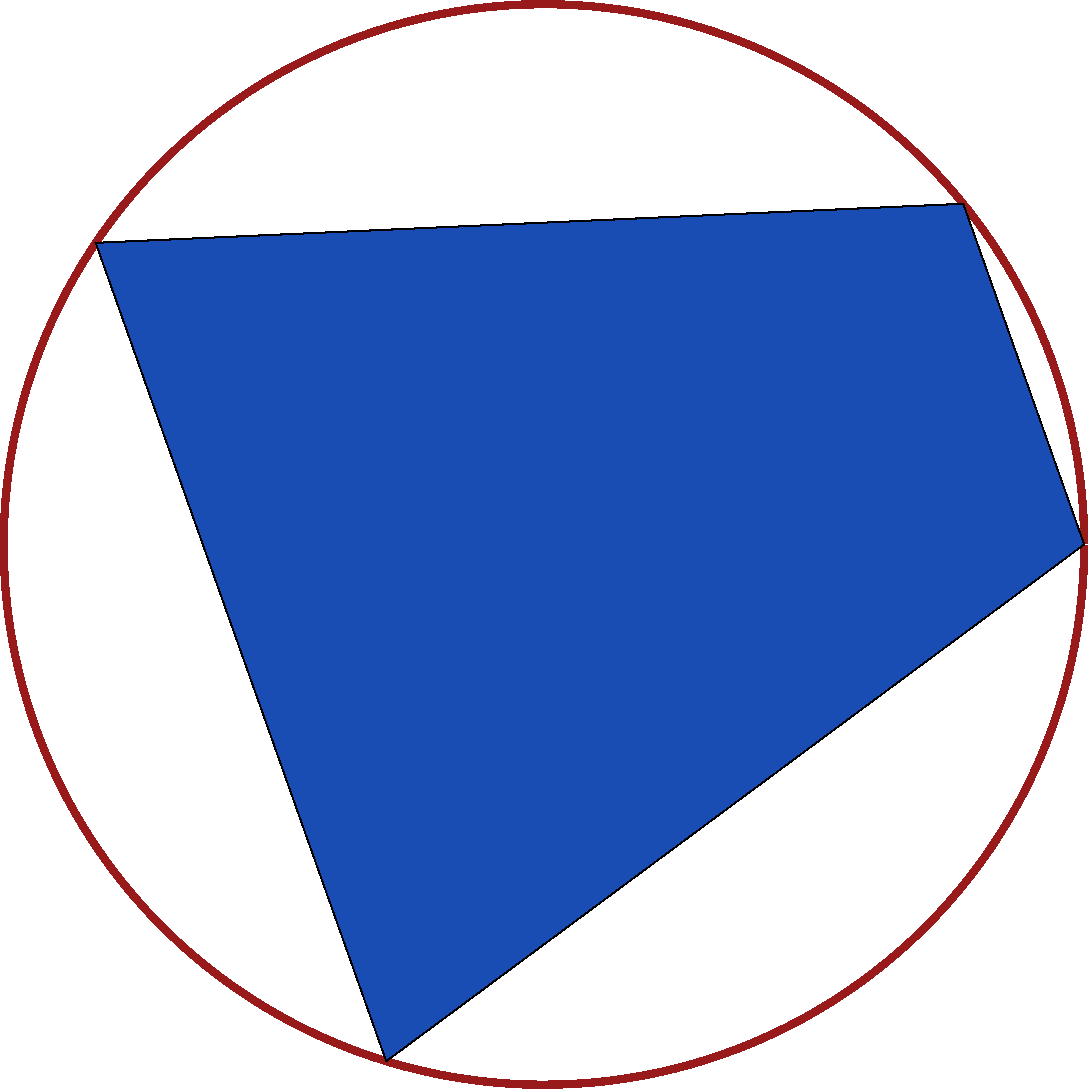} \\
		$A=1$ &
		$A=\frac{3}{2}\sin \frac{2\pi}{3}$ &
		$A=1.75$ \\ \vspace{0.3cm}
		\includegraphics[width=0.3\textwidth]{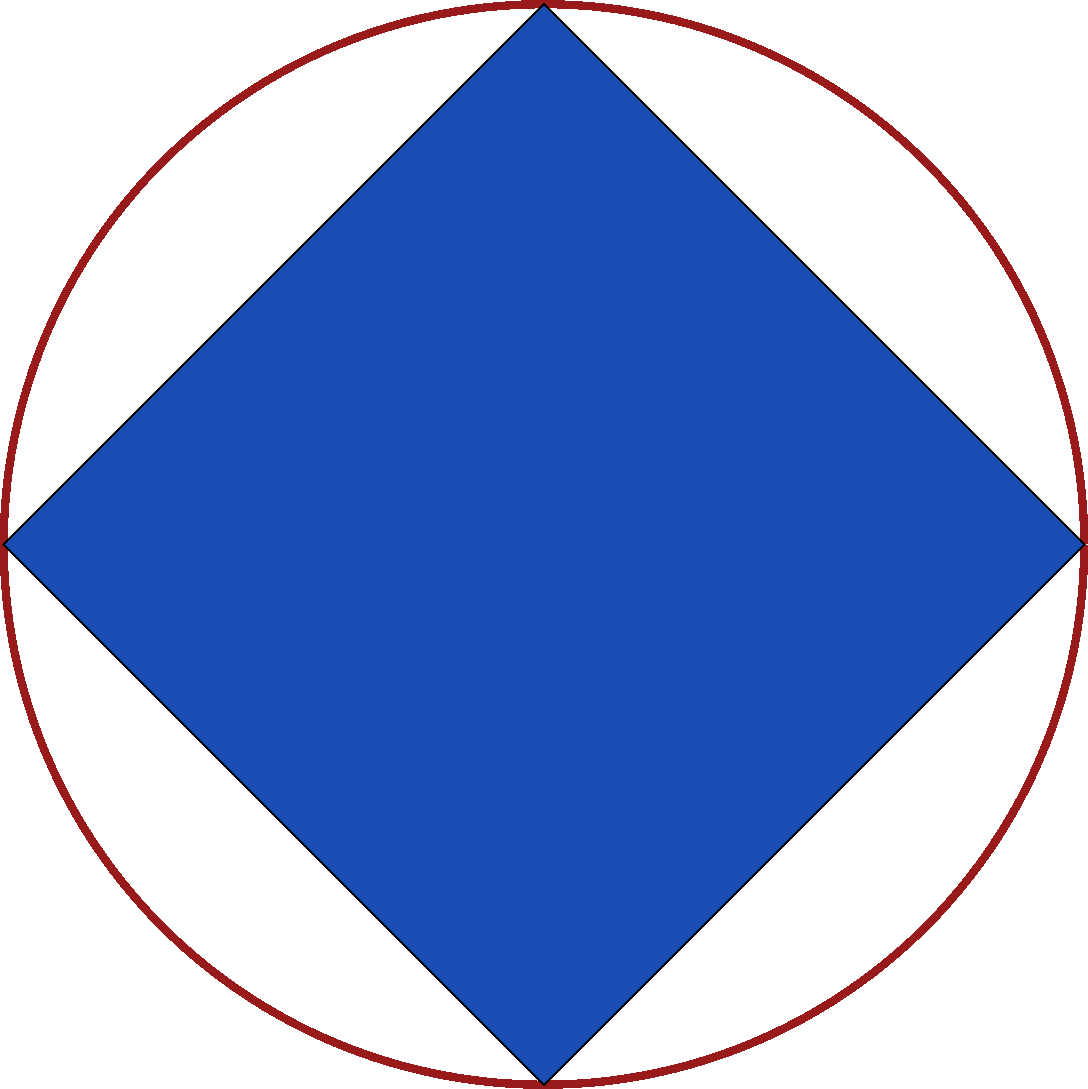}&
		\includegraphics[width=0.3\textwidth]{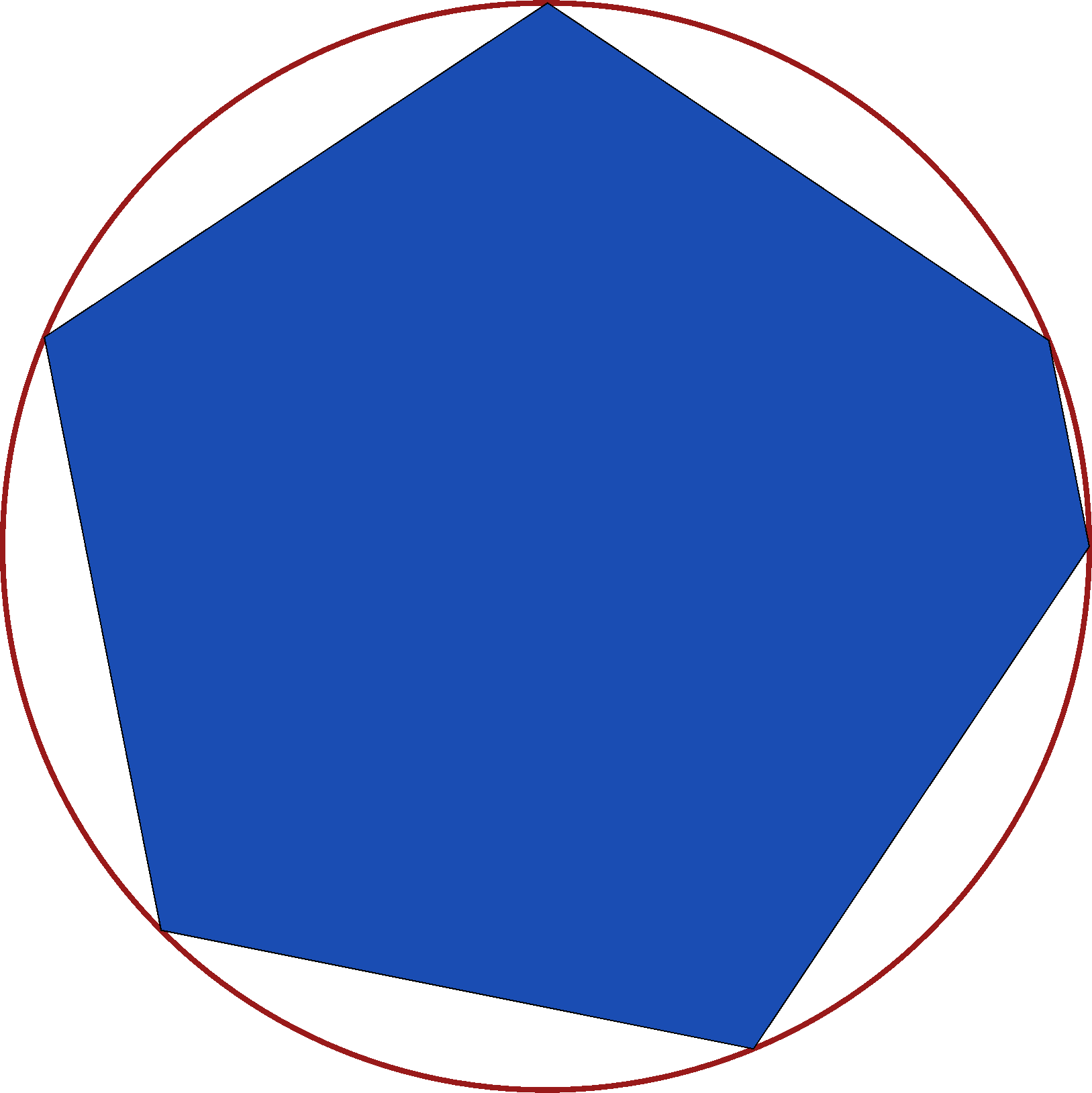}&
		\includegraphics[width=0.3\textwidth]{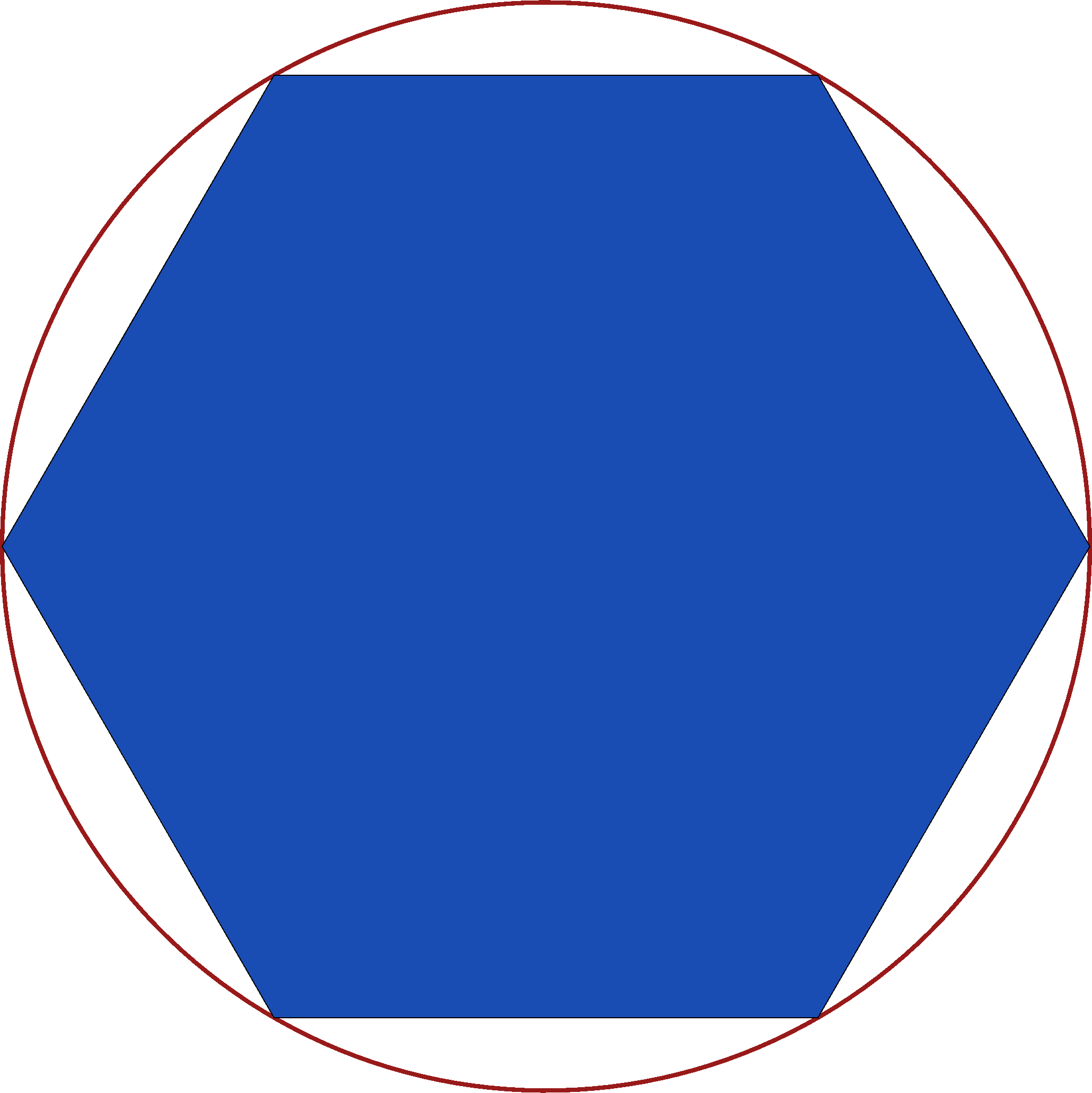} \\
		$A=2$ &
		$A=2.5$ &
		$A=\frac{6}{2}\sin \frac{2\pi}{6}$ \\
		\includegraphics[width=0.3\textwidth]{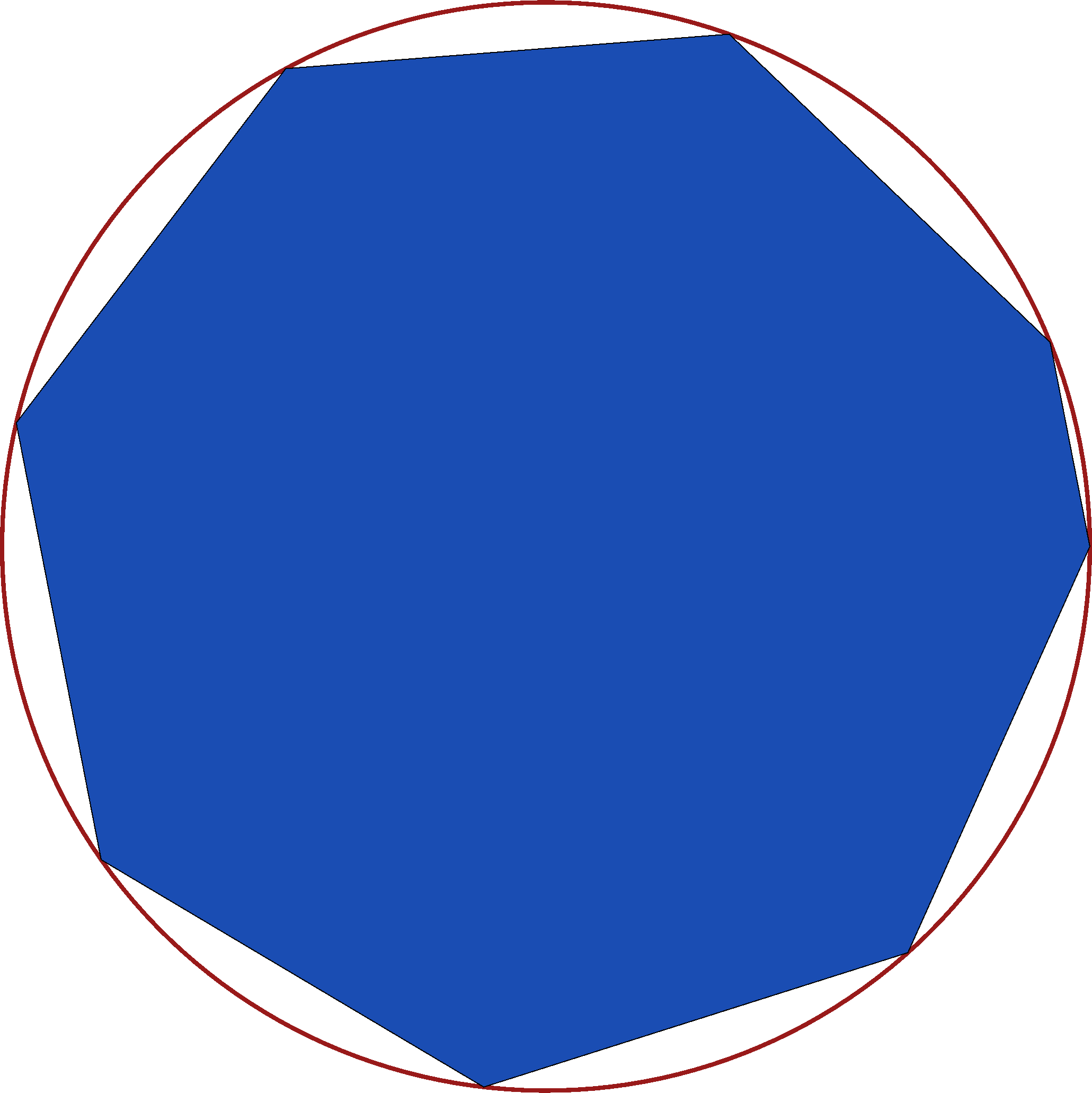}&
		\includegraphics[width=0.3\textwidth]{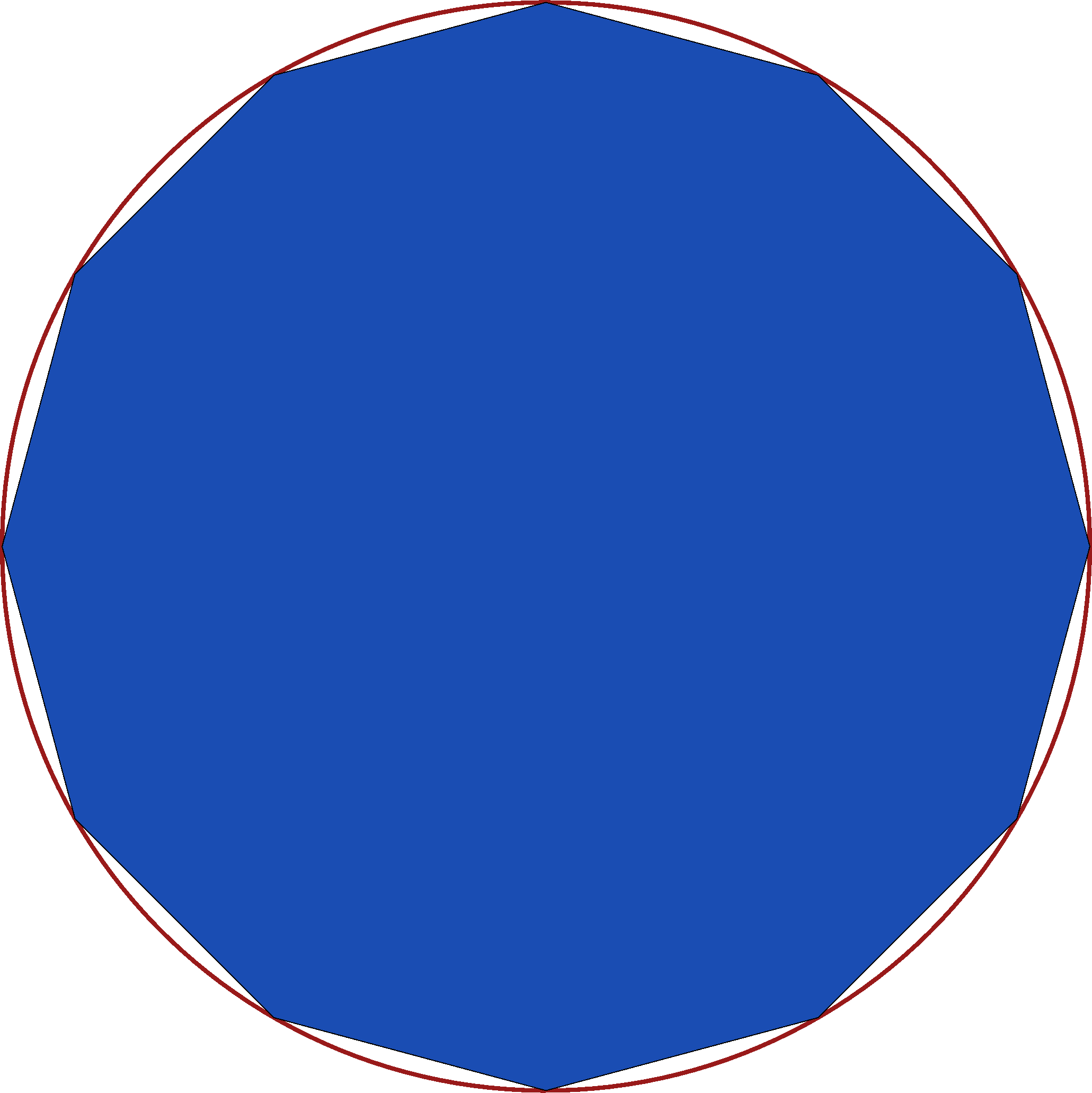}&
		\includegraphics[width=0.3\textwidth]{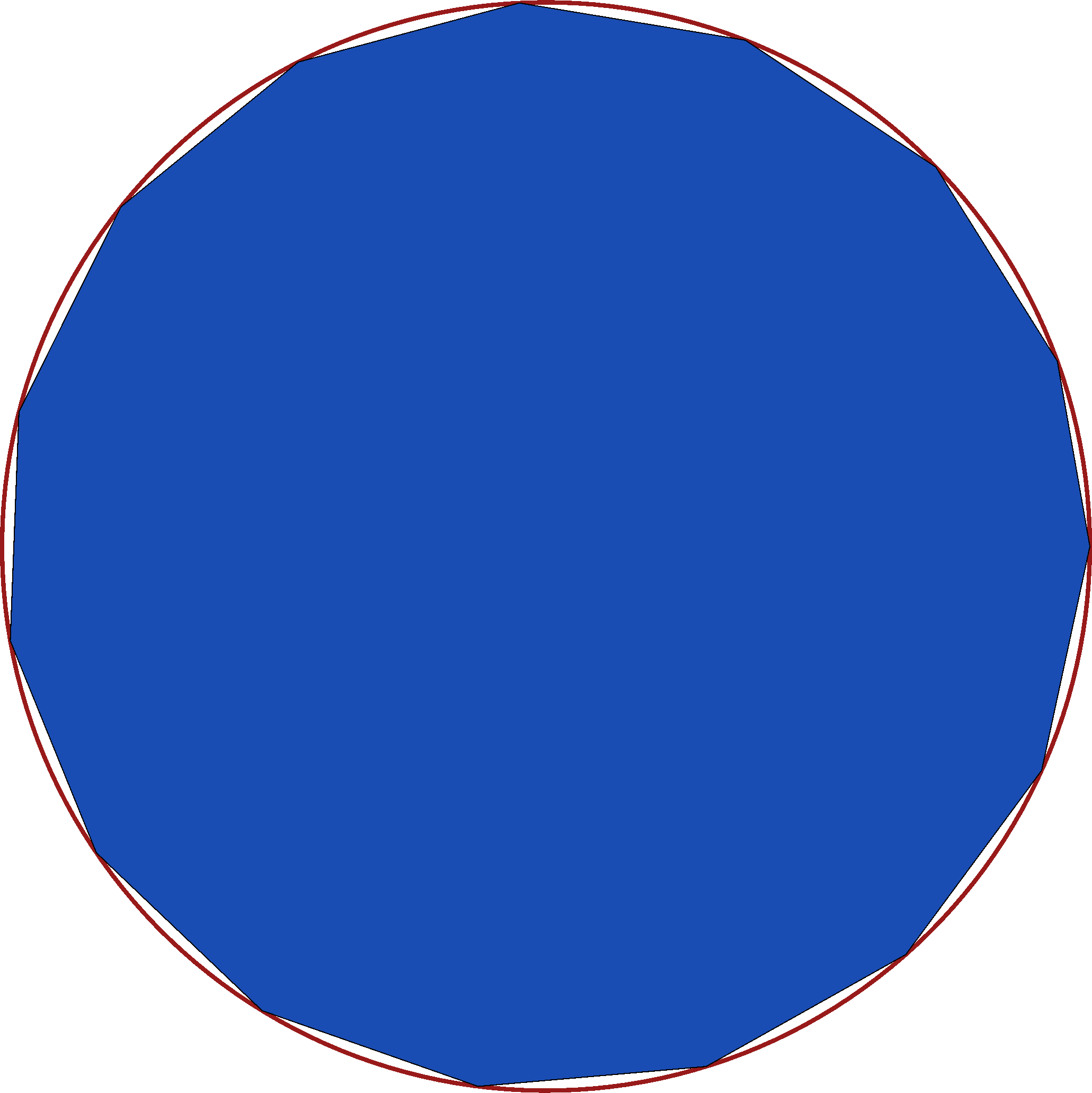} \\
		$A=2.8$ &
		$A=3$ &
		$A=3.05$
	\end{tabular}
\caption{Maximal perimeter convex sets contained in the unit disk $D$ for different area fractions $A \in (0,|D|)$. Changes in the number of edges occur when $A$ is the area of a regular $n$-gon inscribed in $D$.}
\label{fig:solutions-disk} 
\end{figure}

The above steps can be implemented using a simple root finding algorithm for finding zeros in an interval for a monotone function. Examples of optimal polygons for \eqref{eq:max-per} are shown in Figure \ref{fig:solutions-disk} for various values of the area fraction $A \in (0,\pi)$.

\subsection{General convex containers}

In this section a numerical framework is proposed for maximizing the perimeter of a polygon inscribed in a regular convex container $\Omega$ and having a given area. For simplicity, $\Omega$ is parametrized by its radial function $\rho_\Omega : [0,2\pi]\to \Bbb{R}_+$. Of course, it is assumed that $\Omega$ contains the origin. Denote $\bo r(\theta) = (\cos \theta, \sin \theta)$ the radial vector in the direction $\theta$.

Given $n \geq 3$ the number of vertices of the polygon, the algorithm searches $n$ angles $\theta_1\leq \theta_2\leq ... \leq \theta_n \in [0,2\pi]$ such that the polygon with vertices $\bo a_i=\rho_\Omega(\theta_i) \bo r(\theta_i)$, $i=1,...,n$ has a given area $A \in (0,|\Omega|)$ and maximizes the perimeter. Denoting by $\bo o$ the origin, in the implementation the following quantities are computed:
\[ \text{Area}(\Delta \bo a_i \bo o \bo a_{i+1}) = \frac{1}{2} \rho_\Omega(\theta_i) \rho_\Omega(\theta_{i+1}) \sin (\theta_{i+1}-\theta_i)\]
\[ |\bo a_i \bo a_{i+1}| = \sqrt{\rho_\Omega(\theta_i)^2+\rho_\Omega(\theta_{i+1})^2-2\cos (\theta_{i+1}-\theta_i) \rho_\Omega(\theta_i)\rho_\Omega(\theta_{i+1})}. \]
These quantities are enough to compute the area and perimeter of the polygon $\bo a_1,\bo a_2,...,\bo a_n$. Of course, assuming $\rho_\Omega$ is $C^1$ the partial derivatives can also be computed with respect to all variables $(\theta_i)_{i=1}^n$. Thus, gradient based constrained optimization algorithms can be employed to study the problem. In the computation below the Matlab \texttt{fmincon} routine is used with algorithms \texttt{interior-point} and \texttt{sqp}.

Taking $\theta_i$ as variables in the optimization is possible algorithm may lead to difficulties since all inequalities of the form $\theta_i\leq \theta_{i+1}$ need to be imposed. To simplify the problem and the implementation, consider the following variables instead:
\[ \theta_1, \delta_2 = \theta_2-\theta_1, \delta_3 = \theta_3-\theta_2, ..., \delta_{n} = \theta_n-\theta_{n-1}.\]
Of course, the computation of the area and the perimeter is straightforward using these new variables. The monotonicity constraints for the radial angles are transformed into positivity constraints for $\theta_1,\delta_2,..., \delta_n$. Also, the last vertex $\bo a_n$ corresponding to the angle $\theta_n = \theta_1+\delta_2+...+\delta_n$ should not go beyond $\bo a_1$ on $\partial \Omega$, therefore
\begin{equation}\label{eq:constraints-delta} \delta_2+...+\delta_n\leq 2\pi.
\end{equation}
Thus, using the new variables the objective function and the constraint and their derivatives can be computed and the constraint \eqref{eq:constraints-delta} together with $\theta_1\geq 0$, $\delta_i \geq 0$, $i=2,...,n$ should be imposed. To stabilize the algorithm, the upper bound $2\pi$ is also consider for all variables. 

The problem is more complex, with many local minima, when considering general containers. Therefore, multiple runs of the algorithm with different random initializations for the optimization variables are used. The result giving the largest perimeter is shown in the figures below. 

Applying the proposed algorithm to the disk gives precisely the results predicted by the theoretical aspects shown in the previous sections.

\begin{figure}
	\centering 
	\begin{tabular}{ccc}
		\includegraphics[width=0.3\textwidth]{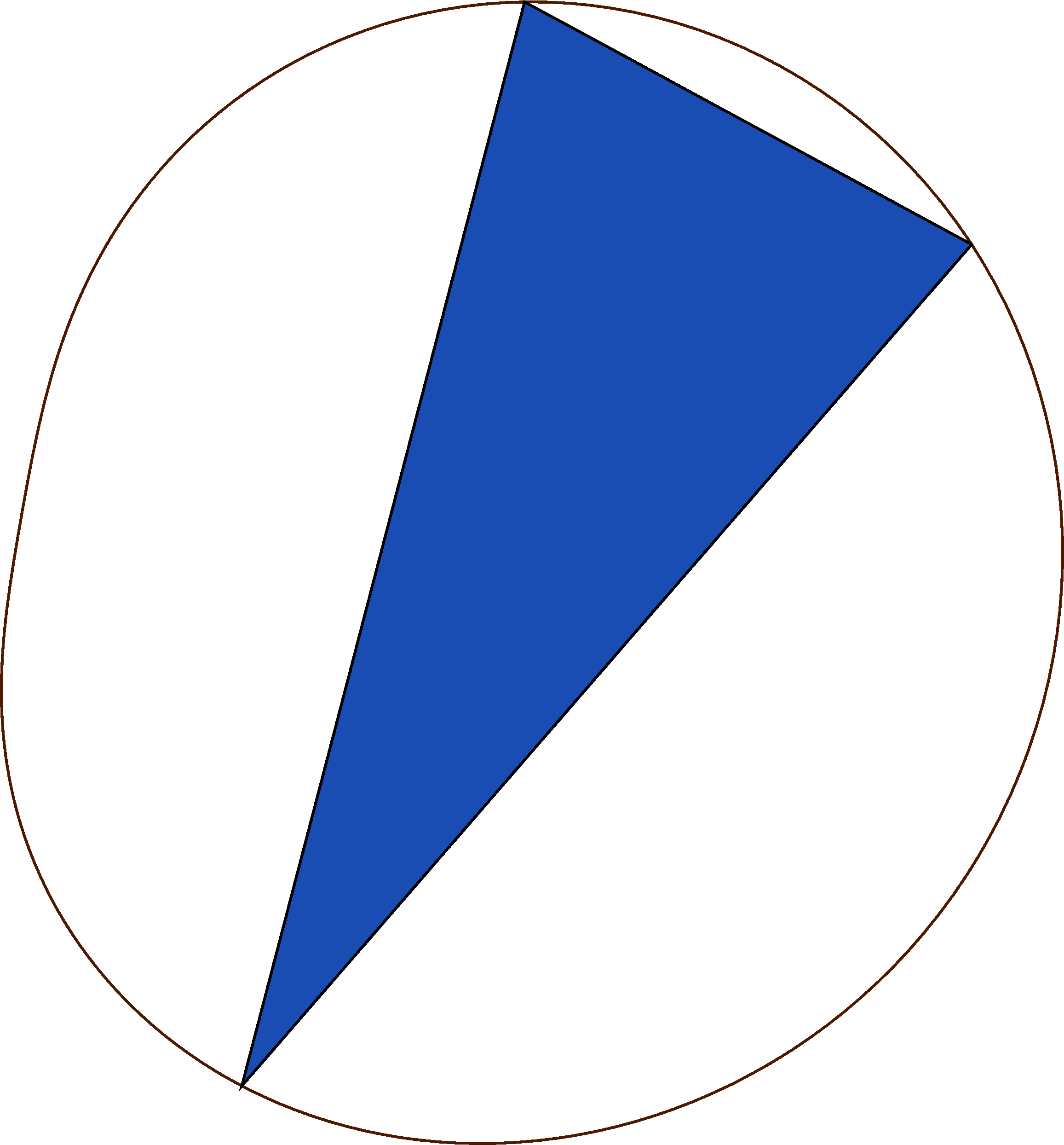}&
		\includegraphics[width=0.3\textwidth]{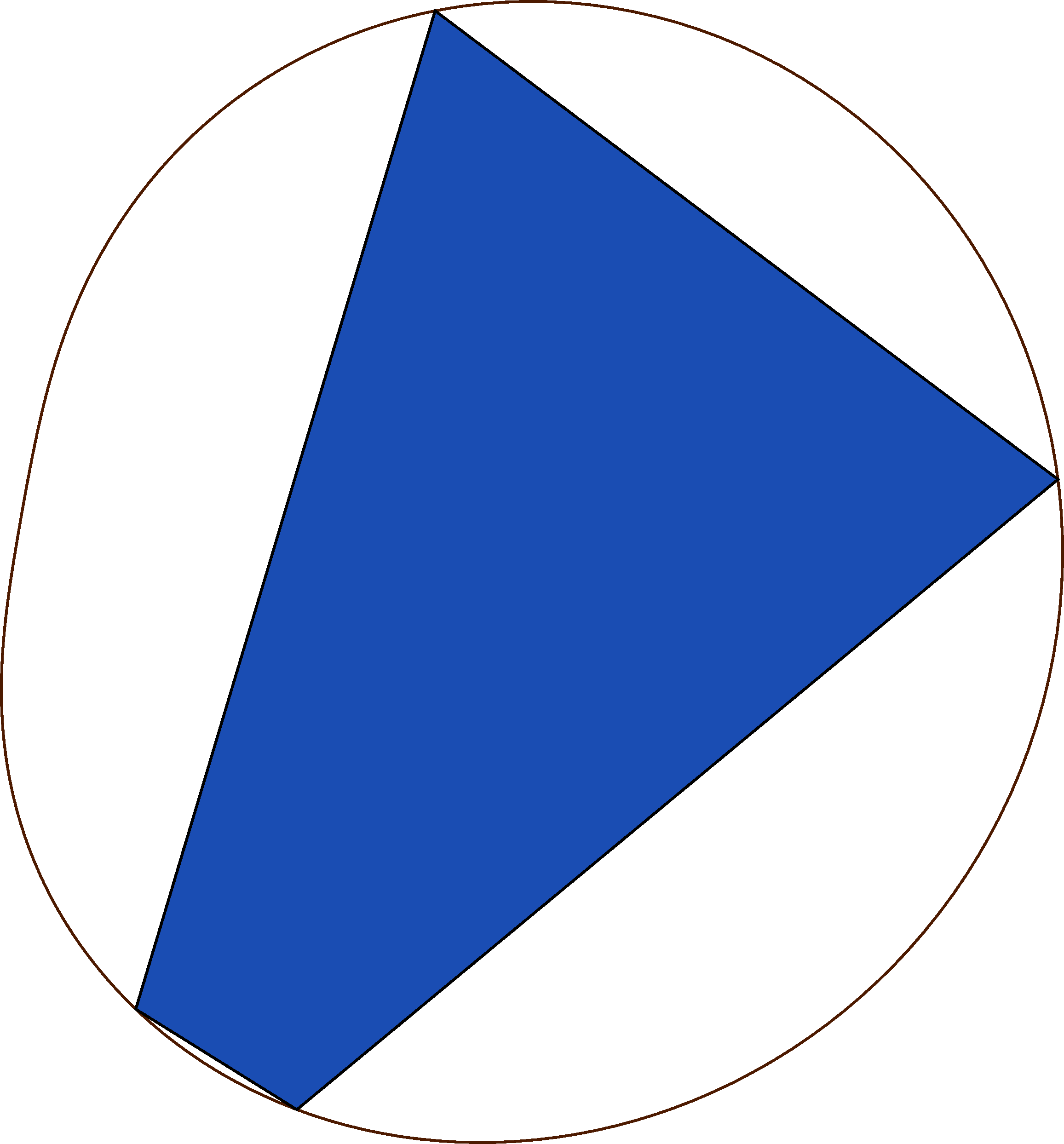}&
		\includegraphics[width=0.3\textwidth]{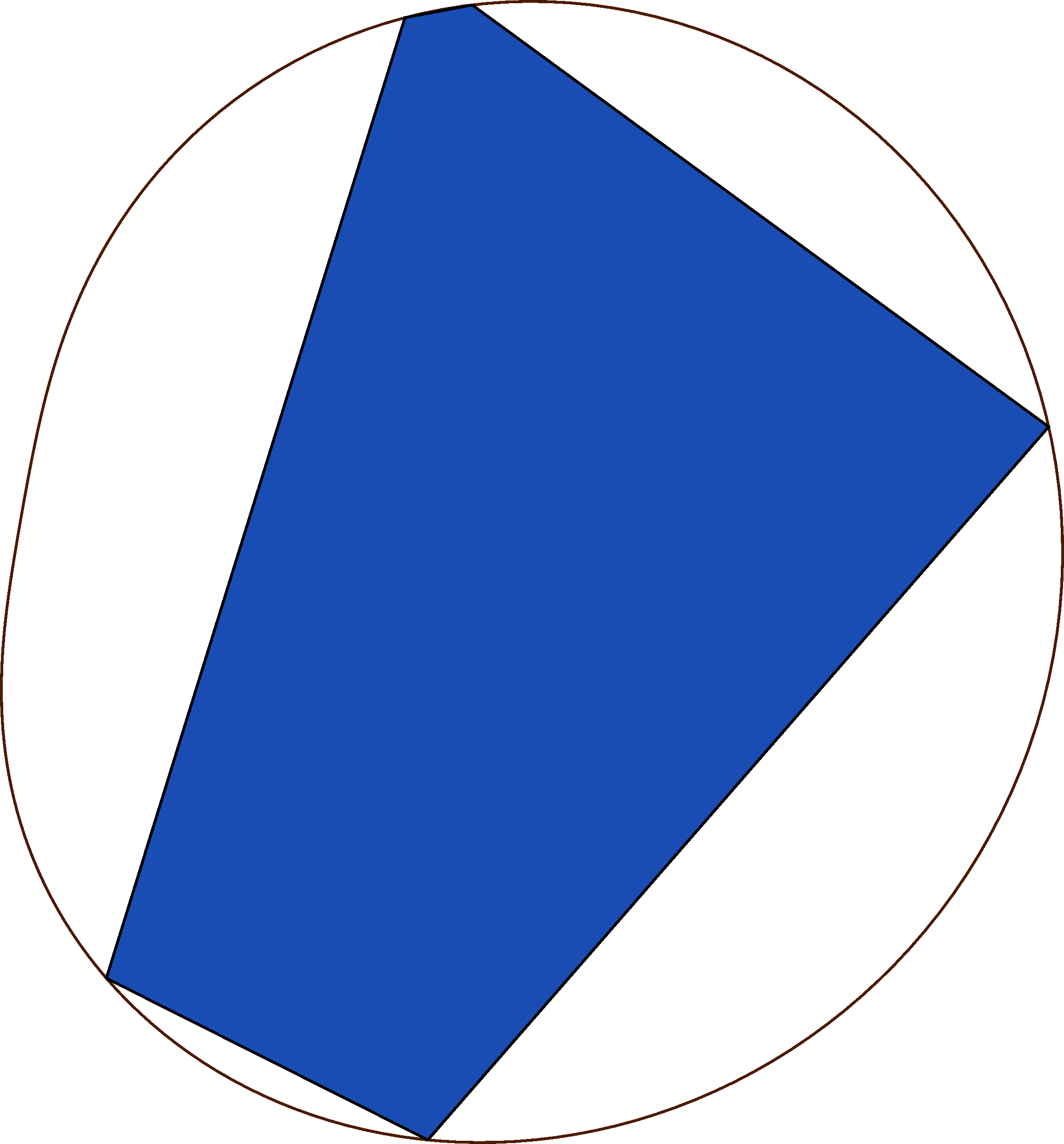} \\
		$A=1$ &
		$A=1.7$ &
		$A=2$ \\ \vspace{0.3cm}
	\end{tabular}
\caption{Numerical results for maximization of the perimeter among polygons with $10,30$ and $50$ sides inscribed in $\Omega_1$ (with radial function given by \eqref{eq:rad_general_1}) with given area fraction $A$.}
\label{fig:general-1}
\end{figure}

Next, consider the shape $\Omega_1$ with radial parametrization
\begin{equation}\label{eq:rad_general_1} \rho_{\Omega_1}(\theta) = 1+0.45\cos \theta+0.04\sin (2\theta).\end{equation}
Optimization results for area fractions $A \in \{1,1.7,2\}$ are shown in Figure \ref{fig:general-1}. The numerical optimization is done in each case among polygons with $10$, $30$ and $50$ vertices. The numerical optimizer is the same polygonal results with $3, 4$ and $5$ vertices, respectively. The optimization algorithm is ran $20$ times with random initializations. This suggests that the result obtained for the disk should extend for general convex containers.

\begin{figure}
	\centering 
	\begin{tabular}{ccc}
		\includegraphics[width=0.3\textwidth]{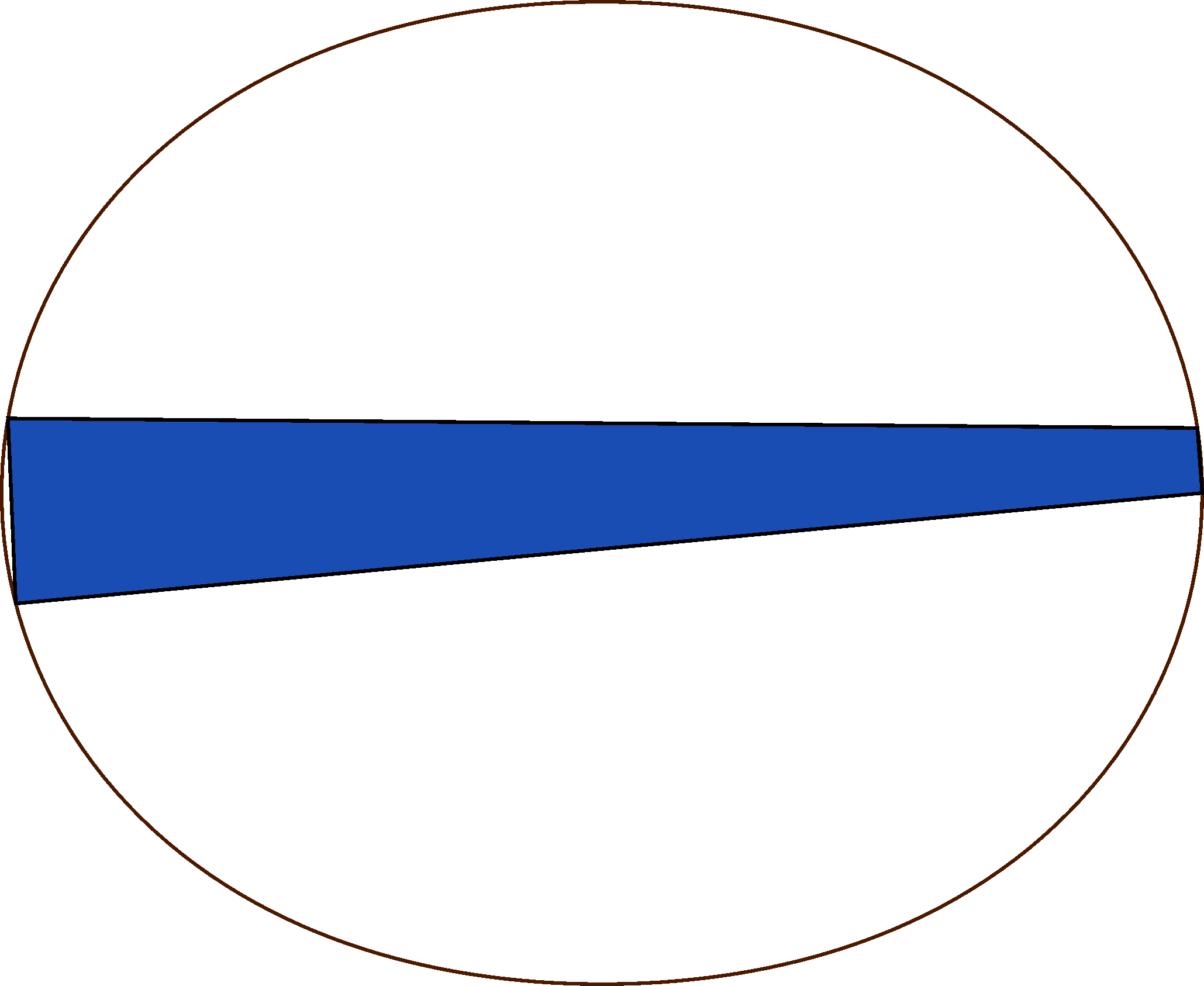}&
		\includegraphics[width=0.3\textwidth]{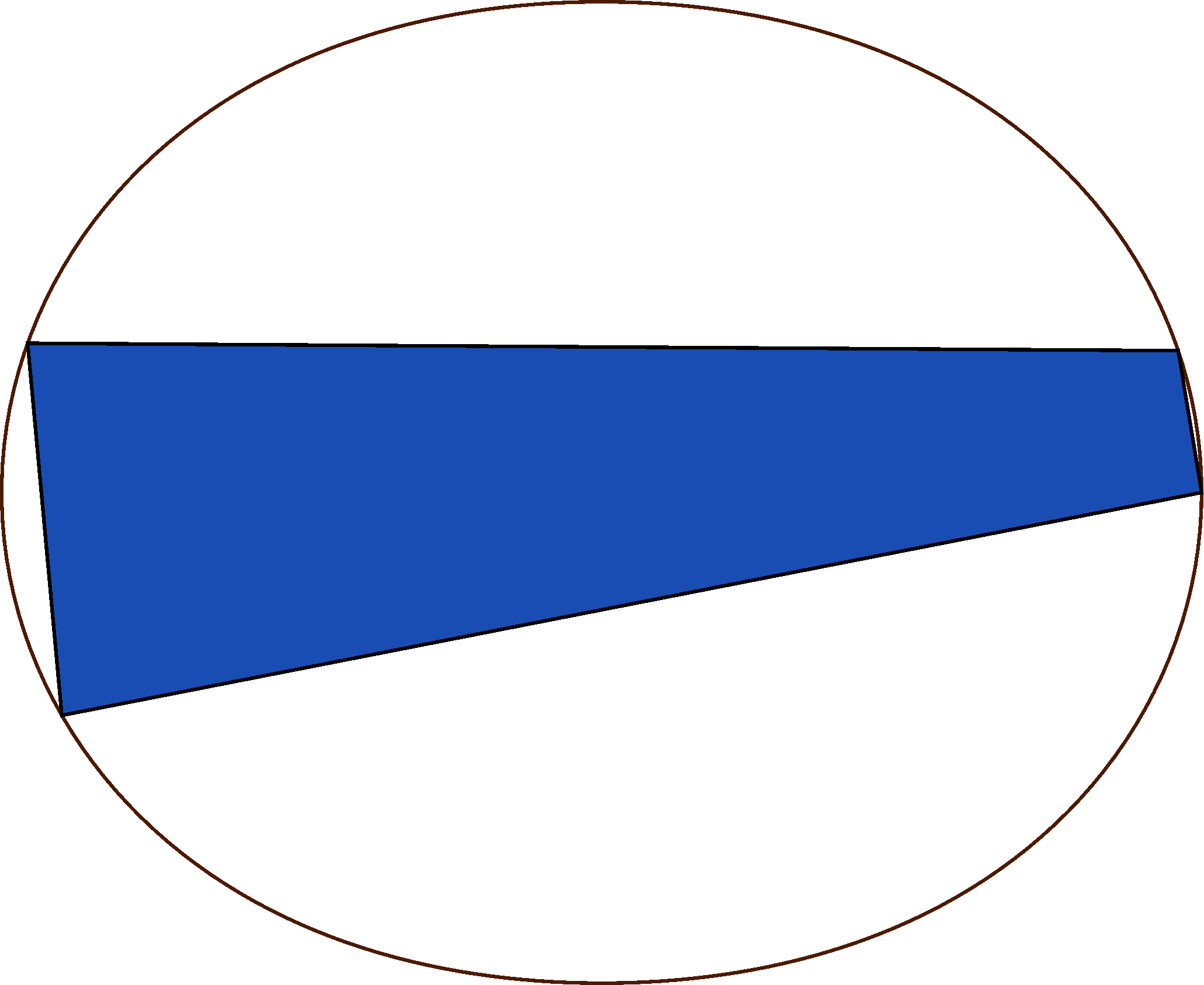}&
		\includegraphics[width=0.3\textwidth]{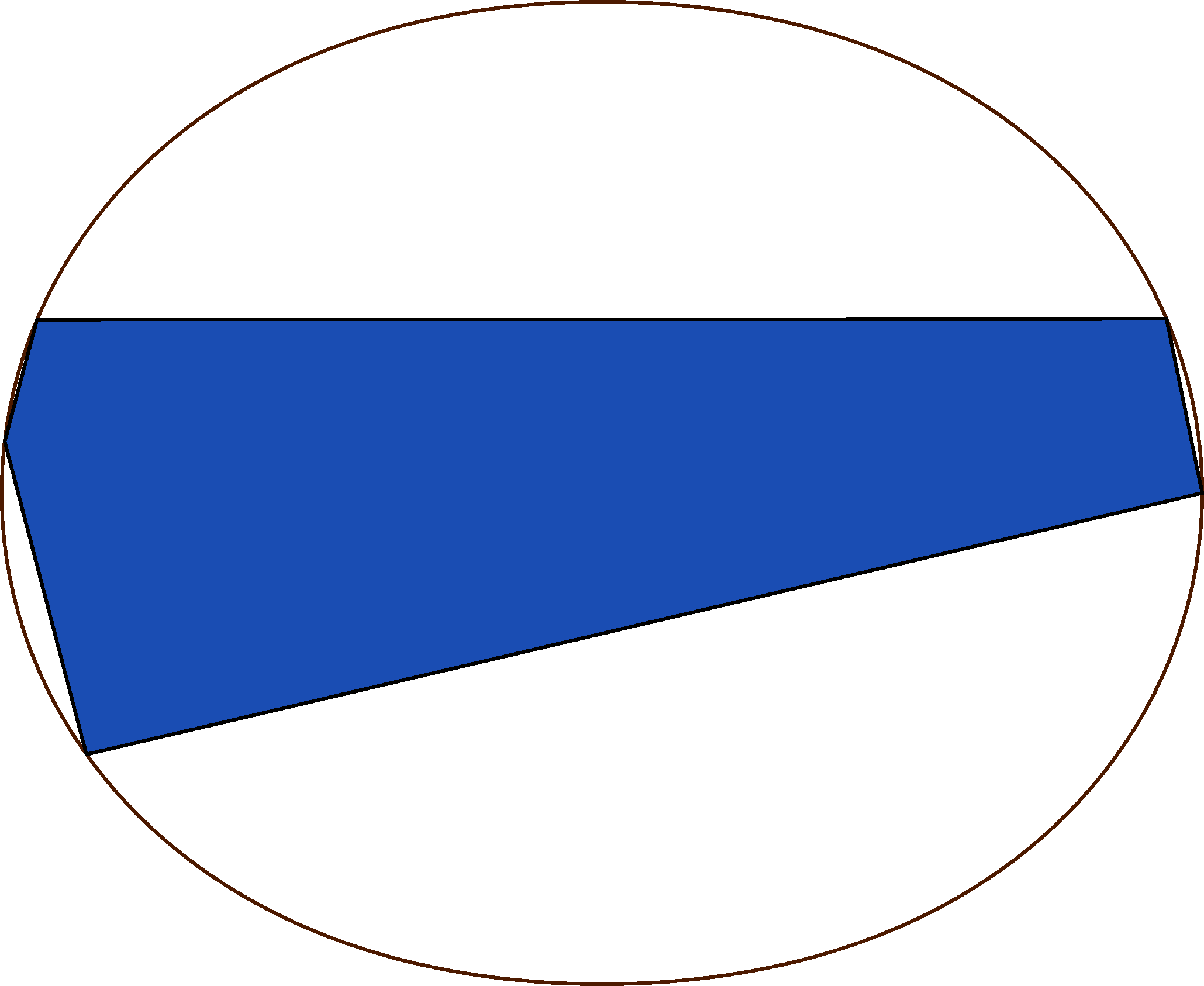} \\
		$A=0.5$ &
		$A=1$ &
		$A=1.2$ \\ \vspace{0.3cm}
	\end{tabular}
	\caption{Numerical results for maximization of the perimeter among polygons with $10,30$ and $50$ sides inscribed in $\Omega_1$ (with radial function given by \eqref{eq:rad_general_2}) with given area fraction $A$.}
	\label{fig:general-2}
\end{figure}

A second set of computations is presented in Figure \ref{fig:general-2} for the shape $\Omega_2$ with radial parametrization given by
\begin{equation}\label{eq:rad_general_2} \rho_{\Omega_1}(\theta) = 1+0.1\cos(2\theta).\end{equation}
and areas $A \in \{0.5,1,1.2\}$. Solutions are again polygonal.

In the general case we have the following observations:
\begin{itemize}[noitemsep]
	\item When $A \to 0$ the optimal shape converges to a diameter of the container $\Omega$. Therefore, for elongated domains it is possible that solutions will have both long edges approaching the diameter and small edges close to the regions of $\partial \Omega$ which are touched by a diameter.
	\item The situation of the disk is very particular: the symmetry allows arbitrary permutations of the sides for the optimizer. The fact that for the disk we have a series of "large" edges and a "small" one is unlikely to generalize to the general case. 
	\item The numerical algorithm proposed for general containers behaves well for small area fractions. For larger area fractions the algorithm faces more and more local minima and results are not as reliable.
\end{itemize}

\section{Conclusions}

In this paper the problem of finding the convex shape contained in a disk, having prescribed area and maximal perimeter is completely solved. The result also appears in \cite{Favard} with an incomplete argument regarding the polygonal character of the minimizer. The optimal shapes are always polygonal with all but one sides equal and completely characterized by the area fraction. The proof is achieved by investigating the analogue problem in the class of polygons, observing that solutions do not change when allowing more than a well characterized number of sides. 

The general case is discussed both from theoretical and numerical points of view. When the containing shape is polygonal, the perimeter maximizer under area constraint is also polygonal. When the containing shape is a general convex set, theoretical aspects become more challenging. A numerical algorithm is proposed and implemented showing that in certain situations it seems that perimeter maximizers under area constraint are indeed polygonal.

\bibliography{./DiscreteIsop}
\bibliographystyle{abbrv}

		\noindent Beniamin \textsc{Bogosel}, Centre de Math\'ematiques Appliqu\'ees, CNRS, École Polytechnique, Institut Polytechnique de Paris, 91120 Palaiseau, France\\
		\nolinkurl{beniamin.bogosel@polytechnique.edu}

\end{document}